\input amstex
\documentstyle{amsppt}
\pagewidth{5.2in}\vsize8.5in\parindent=6mm
\parskip=3pt\baselineskip=14pt\tolerance=10000\hbadness=500
\rightheadtext{the maximal quasiradial Bochner-Riesz operator}
\document
\topmatter
\title
Weak type estimates of the maximal quasiradial Bochner-Riesz
operator on certain Hardy spaces
\endtitle
\author Yong-Cheol Kim \endauthor
\abstract Let $\{A_t\}_{t>0}$ be the dilation group in ${\Bbb
R}^n$ generated by the infinitesimal generator $M$ where
$A_t=\exp(M\log t)$, and let $\varrho\in C^{\infty}({\Bbb
R}^n\setminus\{0\})$ be a $A_t$-homogeneous distance function
defined on ${\Bbb R}^n$. For $f\in {\frak S}({\Bbb R}^n)$, we
define the maximal quasiradial Bochner-Riesz operator ${\frak
M}^{\delta}_{\varrho}$ of index $\delta>0$ by $${\frak
M}^{\delta}_{\varrho} f(x)=\sup_{t>0}\left|{\Cal
F}^{-1}[(1-\varrho/t)_+^{\delta}\hat f\,](x)\right|.$$

If $A_t=t\,I$ and $\{\xi\in {\Bbb R}^n|\,\varrho(\xi)=1\}$ is a
smooth convex hypersurface of finite type, then we prove in an
extremely easy way that ${\frak M}^{\delta}_{\varrho}$ is well
defined on $H^p({\Bbb R}^n)$ when $\delta=n(1/p-1/2)-1/2$ and
$0<p<1$; moreover, it is a bounded operator from $H^p({\Bbb R}^n)$
into $L^{p,\infty}({\Bbb R}^n)$.

If $A_t=t\,I$ and $\varrho\in C^{\infty}({\Bbb R}^n\setminus\{0\})$, we
also prove that ${\frak M}^{\delta}_{\varrho}$ is a bounded
operator from $H^p({\Bbb R}^n)$ into $L^p({\Bbb R}^n)$ when
$\delta>n(1/p-1/2)-1/2$ and $0<p<1$.
\endabstract
\thanks 2000 Mathematics Subject Classification: 42B15, 42B25.
\endthanks
\thanks The author was supported in part by Korea Research Foundation
Proj. No. 2000-003-D00011 and KOSEF Proj. No. 2000-1-10100-001-3.
\endthanks
\address Department of Mathematics Ed., Korea University, Seoul 136-701, Korea
\endaddress
\email ychkim$\@$korea.ac.kr \endemail
\endtopmatter

\define\supp{{\text{\rm supp }}}

\define\inn#1#2{\langle#1,#2\rangle}

\define\lcontr{\rfloor}
\define\lco#1#2{{#1}\lcontr{#2}}
\define\lcoi#1#2{\imath({#1}){#2}}
\define\rco#1#2{{#1}\rcontr{#2}}
\redefine\exp{{\text{\rm exp}}}

\define\ap{\alpha}             
\define\bt{\beta}
\define\gm{\gamma}             
\define\dt{\delta}             

\define\zt{\zeta}
\define\th{\theta}             \define\Th{\Theta}

\define\kp{\kappa}

\define\ld{\lambda}            
\define\sm{\sigma}             \define\Sm{\Sigma}

             \define\Om{\Omega}
\define\vr{\varrho}            \define\iy{\infty}
\define\lt{\left}            \define\rt{\right}
\define\f{\frac}             \define\el{\ell}

\define\fM{{\frak M}}

\define\fR{{\frak R}}
\define\fS{{\frak S}}

\define\fa{{\frak a}}
\define\fb{{\frak b}}

\define\fh{{\frak h}}

\define\BN{{\Bbb N}}

\define\BR{{\Bbb R}}

\define\BZ{{\Bbb Z}}

\define\cB{{\Cal B}}

\define\cD{{\Cal D}}
\define\cE{{\Cal E}}
\define\cF{{\Cal F}}

\define\cH{{\Cal H}}

\define\cN{{\Cal N}}
\define\cO{{\Cal O}}

\define\cQ{{\Cal Q}}

\define\cY{{\Cal Y}}




\define\fhat{{\hat f}}

\define\la{\langle}          \define\ra{\rangle}
         
\define\n{\nabla}            \define\e{\eta}
        \define\fd{\fallingdotseq}
       
     \define\ds{\dsize}
        \define\ls{\lesssim }
\define\gs{\gtrsim}

\subheading{1. Introduction}

Let $\fS(\BR^n)$ be the Schwartz space on $\BR^n$. For
$f\in\fS(\BR^n)$, we denote the Fourier transform of $f$ by $$\cF
[f](x)=\fhat (x)=\int_{\BR^n} e^{-i\la x,\xi\ra} f(\xi)\,d\xi.$$
Then the inverse Fourier transform of $f$ is given by
$$\cF^{-1}[f](x)=\check f (x)=\f{1}{(2\pi)^n}\int_{\BR^n} e^{i\la
x,\xi\ra} f(\xi)\,d\xi.$$ Let $M$ be a real-valued $n\times n$
matrix whose eigenvalues have positive real parts. Then we
consider the dilation group $\{A_t\}_{t>0}$ in $\BR^n$ generated
by the infinitesimal generator $M$, where $A_t=\exp(M\log t)$ for
$t>0$. We introduce $A_t$-homogeneous distance functions $\vr$
defined on $\BR^n$; that is, $\vr :\BR^n\to [0,\iy)$ is a
continuous function satisfying $\vr(A_t\xi)=t\vr(\xi)$ for all
$\xi\in\BR^n$. One can refer to [3] and [11] for its fundamental
properties.

In what follows we shall denote by
$\Sm_{\vr}\fd\{\xi\in\BR^n|\,\vr(\xi)=1\}$ the unit sphere of
$\vr$ and denote by $\BR^n_0=\BR^n\setminus\{0\}$. We use the
polar coordinates; given $x\in\BR^n$, we write $x=r\th$ where
$r=|x|$ and $\th=(\th_1,\th_2,\cdots,\th_n)\in S^{n-1}$. Given
two quantities $A$ and $B$, we write $A\ls B$ or $B\gs A$ if
there is a positive constant $c$ ( possibly depending on the
dimension $n$ and the index $p$ to be given ) such that $A\le c
B$. We also write $A\sim B$ if $A\ls B$ and $B\ls A$.

For $f\in\fS(\BR^n)$, we consider quasiradial Bochner-Riesz means
of index $\dt>0$ defined by $$\fR^{\dt}_{\vr,t}
f(x)=\cF^{-1}[(1-\vr/t)_+^{\dt}\fhat\,](x),$$ and the
corresponding maximal operator $$\fM^{\dt}_{\vr}
f(x)=\sup_{t>0}\lt|\fR^{\dt}_{\vr,t} f(x)\rt|.$$ In the special
case that $\vr(\xi)=|\xi|^2$ and $A_t=t\,I$, Stein, Taibleson,
and Weiss [10] proved that if $0<p<1$, then $\fM^{\dt}_{\vr}$ is
bounded from $H^p(\BR^n)$ into $L^{p,\iy}(\BR^n)$ at the critical
index $\dt=\dt(p)\fd n(1/p-1/2)-1/2$ where $H^p(\BR^n)$ is the
standard real Hardy space defined in Stein [9] and
$L^{p,\iy}(\BR^n)$ is one of the Lorentz spaces (which is called
weak-$L^p$ space) defined in Stein and Weiss [12] and furthemore
Stein obtained the exceptional result that there is $f\in H^1
(\BR^n)$ such that a.e. convergence of the Bochner-Riesz means
fails for $p=1$ and $\dt(1)=(n-1)/2$.

In our first result we shall assume that $\vr\in
C^{\iy}(\BR^n_0)$, $A_t=t\,I$ and $\Sm_{\vr}$ is a smooth convex
hypersurface of $\BR^n$ which is of finite type, i.e. every
tangent line makes finite order of contact with $\Sm_{\vr}$. We
say that $\Sm_{\vr}$ is of finite type $k\ge 2$ if $k$ is the
maximal order of contact on $\Sm_{\vr}$.

\proclaim{Theorem 1.1} Suppose that $A_t=t\,I$, $\vr\in
C^{\iy}(\BR^n_0)$ is a $A_t$-homogeneous distance function defined
on $\BR^n$, and $\Sm_{\vr}$ is a smooth convex hypersurface of
finite type. Then $\fM^{\dt(p)}_{\vr}$ is well defined on
$H^p(\BR^n)$ when $0<p<1$; moreove, $\fM^{\dt(p)}_{\vr}$ is a
bounded operator from $H^p(\BR^n)$ into $L^{p,\iy}(\BR^n)$. That
is, there is a constant $C=C(n,p,\Sm_{\vr})>0$ such that for any
$f\in H^p(\BR^n)$, $$\lt|\{x\in\BR^n|\,\fM^{\dt(p)}_{\vr}
f(x)>\ld\}\rt|\le \f{C}{\ld^p}\,\|f\|_{H^p(\BR^n)}^p,\,\,\ld>0,$$
where $|E|$ denotes the Lebesgue measure of the set
$E\subset\BR^n$.
\endproclaim

\noindent{\it Remark.} As a matter of fact, we prove this result
under more general surface condition than the finite type
condition on $\Sm_{\vr}$, which is to be called a spherically
integrable condition of order $<1$ in Section 3.

Our second result is to obtain that if $\dt>n(1/p-1/2)-1/2$ and
$0<p<1$ then $\fM^{\dt}_{\vr}$ admits $(H^p,L^p)$-estimate under
no surface condition on $\Sm_{\vr}$.

\proclaim{Theorem 1.2} Suppose that $A_t=t\,I$ and $\vr\in
C^{\iy}(\BR^n_0)$ is a $A_t$-homogeneous distance function defined
on $\BR^n$. If $\dt>\dt(p)$ for $0<p<1$, then $\fM^{\dt}_{\vr}$ is
a bounded operator from $H^p(\BR^n)$ into $L^p(\BR^n)$; that is,
there is a constant $C=C(n,p)>0$ such that for any $f\in
H^p(\BR^n)$, $$\|\fM^{\dt}_{\vr} f\|_{L^p(\BR^n)}\le
C\,\|f\|_{H^p(\BR^n)},$$ provided that $\dt>n(1/p-1/2)-1/2$ and
$0<p<1$.
\endproclaim

\noindent{\it Remark.} This problem is still left open on the
critical index $\dt=n(1/p-1/2)-1/2$ and $0<p<1$.

\subheading{2. $(H^p,L^p)$-estimate for the case that $\vr\in
C^{\iy}(\BR^n_0)$ and $\dt>\dt(p)$}

We shall employ a decomposition of the Bochner-Riesz multiplier
$(1-\vr)^{\dt}_+$ as in A. C\'ordoba [2].  Let $\phi\in
C^{\iy}_0(1/2,2)$ satisfy $\sum_{k\in\BZ} \phi(2^k t)=1$ for all
$t>0$. For $k\in\BN$, let
$\Phi^{\dt}_k=\phi(2^{k+1}(1-\vr))(1-\vr)_+^{\dt}$ and
$\Phi^{\dt}_0=(1-\vr)^{\dt}_+-\sum_{k\in\BN}\Phi^{\dt}_k$. For
each $k\in\BZ$, we now introduce a partition of unity $\Xi_{k\el},
\el=1,2,\cdots, N_k$, on the unit sphere $\Sm_{\vr}$ which we
extend to $\BR^n$ by way of $\Pi_{k\el}(A_t\zt)=\Xi_{k\el}(\zt),
t>0, \zt\in\Sm_{\vr}$, and which satisfies the following
properties; there are a finite number of points
$\zt_{k1},\zt_{k2},\cdots,\zt_{k N_k}\in\Sm_{\vr}$ such that for
$\el=1,2,\cdots,N_k$,

(i) $\sum_{\el=1}^{N_k}\Pi_{k\el}(\zt)\equiv 1$ for all
$\zt\in\Sm_{\vr}$,

(ii) $\Xi_{k\el}(\zt)=1$ for all $\zt\in \Sm_{\vr}\cap
 B(\zt_{k\el};2^{-k/2})$,

(iii) $\Xi_{k\el}$ is supported in $\Sm_{\vr}\cap
 B(\zt_{k\el};c_1\,2^{-k/2})$,

(iv) $\lt|\cD^{\ap}\Pi_{k\el}(\xi)\rt|\le c_2\,2^{|\ap|k/2}$ for
any multiindex $\ap$, if $1/2\le\vr(\xi)\le 2$,

(v) $N_k\le c_3\,2^{(n-1)k/2}$ for fixed $k$,

\noindent where $B(\zt_0;s)$ denotes the ball in $\BR^n$ with
center $\zt_0\in\Sm_{\vr}$ and radius $s>0$ and the positive
constants $c_1, c_2, c_3$ do not depend upon $k$. For each
$k\in\BZ$, let $\cH^{\dt}_{\vr
k\el}=\cF^{-1}[\Phi^{\dt}_k\Pi_{k\el}]$ and
$\cH_0=\cF^{-1}[\Phi^{\dt}_0]$.

Next we invoke a simple observation used in [8] to obtain decay
estimate for kernels $\cH_{k\el}$, $\cH_0$ corresponding to the
decomposition of the Bochner-Riesz multiplier defined in the
above. Without loss of generality, we can assume that $\vr\in
C^{\iy}(\BR^n)$ because we can replace $\vr$ by $\vr^N$ for
sufficiently large $N>0$ by a subordination argument in [3]. Then
we easily see that the kernel $\cH_0$ has a nice decay, and so
its corresponding maximal operator admits
$(H^p,L^{p,\iy})$-estimate for the critical index
$\dt(p)=n(1/p-1/2)-1/2$ and $0<p<1$ as in that of Stein,
Taibleson, and Weiss [10]. Thus we concentrate upon obtaining the
decay estimate for the kernels $\cH^{\dt}_{\vr k\el}$.

\proclaim{Lemma 2.1.} For fixed $k\in\BN$ and for
$\el=1,2,\cdots,N_k$, let $T_{\zt_{k\el}}(\Sm_{\vr})$ be the
tangent space of $\Sm_{\vr}$ at $\zt_{k\el}\in\Sm_{\vr}$,
$\{e^j_{k\el}\}_{j=1}^{n-1}$ be an orthonormal basis of
$T_{\zt_{k\el}}(\Sm_{\vr})$, and $e_{k\el}^0$ be the outer unit
normal vector to $\Sm_{\vr}$ at $\zt_{k\el}\in\Sm_{\vr}$. Then we
have the following estimate $$\lt|\cH^{\dt}_{\vr k\el}(x)\rt|\le
\f{ C_N\,2^{-k(\dt+1+(n-1)/2)}}{(1+2^{-k}|\la x,e^0_{k\el}\ra|)^N
\prod_{j=1}^{n-1} (1+2^{-k/2}|\la x,e_{k\el}^j\ra|)^N}$$ for any
$N\in\BN$.
\endproclaim

\noindent{\it Proof.} We need the following simple observation:

Let $\vr\in C^N(\BR^n)$ and $F\in C^N(\BR_+)$. For $e\in S^{n-1}$,
let $\cD_e f$ be the directional derivative $\la e,\n f\ra$. Then
one can have the formula ( see [8] )
$$\cD_e^N (F\circ\vr)=\sum_{\nu=1}^N
F^{(\nu)}\circ\vr\sum_{\bt\in\cY_{\nu}^N}\sum_{m=1}^{\nu}
c_{N,\bt_m}\cD^{\bt}_e\vr \tag{2.1} $$ where
$\cY_{\nu}^N=\{\bt|\,\sum_{m=1}^{\nu}\bt_m=N,\,\text{at least
$\nu-\f{N}{2}$ of the numbers $\bt_m$ are equal to $1$}\}$,
$\bt=(\bt_1,\cdots,\bt_{\nu})$ is a multiindex, and
$c_{N,\bt_m}$'s are some constants. For $k\in\BN$, let
$F_k(t)=\phi(2^{k+1}(1-t))(1-t)^{\dt}_+$. Then it follows from
simple computation that
$$F_k^{(\nu)}(t)=(-1)^{\nu}\sum_{i=0}^{\nu} C(\nu,i) C(\dt,\nu-i)
2^{i(k+1)}\phi^{(i)}(2^{k+1}(1-t))(1-t)^{\dt-\nu+i}\tag{2.2}$$
where $C(\nu,i)=\nu(\nu-1)(\nu-2)\cdots (\nu-i+1)$ for positive
integers $\nu$, $i$, and $C(\nu,0)=1$. For fixed $k,\el$, we have
the estimate $$\lt\|\cD^N_{e^0_{k\el}}[\Phi^{\dt}_k
\Pi_{k\el}]\rt\|_{L^1}\le c\, 2^{-k(\f{n+1}{2})} 2^{-k\dt} 2^{k
N} \tag{2.3} $$ for any $N\in\BN$. Since we have the better
estimate $|\cD_{e^j_{k\el}}\vr|\le c\, 2^{-k/2}$ on the support
of $\cF[\cH_{\vr k\el}^{\dt}]$ for fixed $j, k, \el$, it follows
from (2.1) and Taylor's theorem that
$$\lt\|\cD^N_{e^j_{k\el}}[\Phi^{\dt}_k \Pi_{k\el}]\rt\|_{L^1}\le
c\,2^{-k(\f{n+1}{2})} 2^{-k\dt} 2^{k N/2} \tag{2.4}$$ for any
$N\in\BN$. Using the integration by parts, it follows from (2.3)
and (2.4) that $$ \lt|\cH^{\dt}_{\vr k\el}(x)\rt|\le  \f{
C_N\,2^{-(\dt+1+(n-1)/2)k}}{(1+2^{-k}|\la x,e^0_{k\el}\ra|)^N
\prod_{j=1}^{n-1} (1+2^{-k/2}|\la x,e_{k\el}^j\ra|)^N}
\tag{2.5}$$ for any $N\in\BN$. \qed

We now introduce the real Hardy space $H^p(\BR^n)$ defined in
terms of atomic decompositions along the pattern of Stein [9].
For $0<p\le 1$, a function $\fa\in L^{\iy}(\BR^n)$ is called a
$(p,\mu)$-atom centered at $x_0\in\BR^n$ if it satisfies

(i) there is a ball $B(x_0;s)$ with $\supp(\fa)\subset B(x_0;s)$,

(ii) $\|\fa\|_{L^{\iy}}\le |B(x_0;s)|^{-1/p}$, and

(iii) $\ds\int_{\BR^n} \fa(x) x^{\ap}\,dx=0$ for $|\ap|\le\mu$,

\noindent where $\ap=(\ap_1,\ap_2,\cdots,\ap_n)$ is an $n$-tuple
of nonnegative integers and $|\ap|=\ap_1 +\ap_2 +\cdots +\ap_n$.
If $f=\sum_{k=1}^{\iy} c_k \fa_k$ where the $\fa_k$'s are
$(p,\mu)$-atoms and $\{c_k\}\in\el^p$, then $f\in H^p(\BR^n)$ and
$\|f\|^p_{H^p}\lesssim\sum_k |c_k|^p$ and the converse inequality
also holds. Here we note that if $\dt>n(1/p-1/2)-1/2$ then
$\mu=n(1/p'-1)$ is enough for our oncoming estimates where $p'<p$
is a positive number satisfying $\dt=n(1/p'-1/2)-1/2$.

For $f\in\fS(\BR^n)$, $\dt>0$, $k\in\BN$, and
$\el=1,2,\cdots,N_k$, let $$\fM^{\dt}_{\vr k\el}
f(x)=\sup_{t>0}\lt|\cH^{\dt,t}_{\vr k\el}*f(x)\rt|$$ where
$\cH^{\dt,t}_{\vr k\el}(x)=t^n\, \cH^{\dt}_{\vr k\el}(A_t^* x)$,
and let $\fM^{\dt}_{\vr k} f(x)=\sum_{\el=1}^{N_k} \fM^{\dt}_{\vr
k\el} f(x)$.

\proclaim{Lemma 2.2} If $\dt>n(1/p-1/2)-1/2$ for $0<p<1$, let a
positive number $p'<p$ be chosen so that $\dt=n(1/p'-1/2)-1/2$.
For fixed $k\in\BN$ and for $\el=1,2,\cdots,N_k$, let
$T_{\zt_{k\el}}(\Sm_{\vr})$ be the tangent space of $\Sm_{\vr}$
at $\zt_{k\el}\in\Sm_{\vr}$, $\{e^j_{k\el}\}_{j=1}^{n-1}$ be an
orthonormal basis of $T_{\zt_{k\el}}(\Sm_{\vr})$, and
$e_{k\el}^0$ be the outer unit normal vector to $\Sm_{\vr}$ at
$\zt_{k\el}\in\Sm_{\vr}$. Then we have the following estimate
$$\lt|\cH^{\dt}_{\vr k\el}(x)\rt|+\lt|\n\cH^{\dt}_{\vr
k\el}(x)\rt|\le \f{ C_p\,2^{-k(\f{n-1}{2 p'})}}{
\prod_{j=0}^{n-1} (1+|\la x,e_{k\el}^j\ra|)^{1/p'}}\fd
C_{p}\,2^{-k(\f{n-1}{2 p'})} P_{k\el}(x).$$
\endproclaim

\noindent{\it Proof.} This can easily be obtained by choosing
$\dt=n(1/p'-1/2)-1/2$ and $N=1/p'$ in Lemma 2.1. We also observe
that $\n\cH^{\dt}_{\vr k\el}=\varphi * \cH^{\dt}_{\vr k\el}$ for
some $\varphi\in\fS(\BR^n)$.  \qed

\proclaim{Lemma 2.3} If $\dt>n(1/p-1/2)-1/2$ for $0<p<1$, let a
positive number $p'<p$ be chosen so that $\dt=n(1/p'-1/2)-1/2$.
Suppose that $\fa$ is a $(p,n(1/p'-1))$-atom on $\BR^n$ which is
supported in the ball $B(x_0;s)$ with center $x_0\in\BR^n$ and
radius $s>0$. Then there is a constant $C=C(n,p)>0$ such that

$(a)$ $\ds\lt|\fM^{\dt}_{\vr k\el}\fa(x)\rt|\le
C\,s^{-n/p}\,2^{-k(\f{n-1}{2 p'})}
P_{k\el}\lt(\f{x-x_0}{s}\rt)\,\text { for any $x\in
{B(x_0;2s)}^c$,}$

$(b)$ $\ds\lt\|(\fM^{\dt}_{\vr
k\el}\fa)\chi_{{B(x_0;2s)}^c}\rt\|_{L^p}\le C\,2^{-k(\f{n-1}{2
p'})}$,

\noindent where $P_{k\el}(x)$ is the function given in Lemma 2.2.
\endproclaim

\noindent{\it Proof.} (a) We first assume that $\fa$ is a
$(p,n(1/p'-1))$-atom which is supported in the unit ball $B(0;1)$
centered at the origin and let $N\in\BN$ be an integer satisfying
$N<n(1/p'-1)\le N+1$, i.e. $n/(n+N+1)\le p'<n/(n+N)$. If $x\in
{B(0;2)}^c$ and $t>1$, then it easily follows from Lemma 2.2 that
$$\lt|\cH^{\dt,t}_{\vr k\el} *\fa(x)\rt|\le
C\,t^{n(1-1/p')}\,2^{-k(\f{n-1}{2p'})} P_{k\el}(x).$$ Since
$n(1-1/p')<0$, we have that $$\sup_{t>1}\lt|\cH^{\dt,t}_{\vr
k\el} *\fa(x)\rt|\le C\,2^{-k(\f{n-1}{2p'})}
P_{k\el}(x).\tag{2.6}$$

If $x\in {B(0;2)}^c$ and $0<t\le 1$, let $\cQ_{t,x}(y)$ be the
$N$-th order Taylor polynomial of the function
$y\mapsto\cH^{\dt(p)}_{\vr k\el}(A_t^* (x-y))$ expanded near the
origin. Using the moment conditions on the atom $\fa$ and Taylor's
theorem, we obtain the estimate $$ \split \lt|\fM^{\dt,t}_{\vr
k\el} *\fa(x)\rt|&=t^n\lt|\int_{\BR^n} [\cH^{\dt}_{\vr
k\el}(A_t^* (x-y))-\cQ_{t,x}(y)]\fa(y)\,dy\rt| \\
&\ls t^{n+(N+1)}\int_0^1\int_{B(0;1)}\lt|\n^{N+1}\cH^{\dt}_{\vr
k\el} (A_t^*(x-\tau y))\rt|\,dy\,d\tau \\
&\ls t^{n+(N+1)-n/p'}\,2^{-k(\f{n-1}{2 p'})} P_{k\el}(x)\endsplit
$$ because $n+(N+1)-n/p'\ge 0$. Thus we have that $$\sup_{0<t\le
1}\lt|\cH^{\dt,t}_{\vr k\el} *\fa(x)\rt|\ls 2^{-k(\f{n-1}{2 p'})}
P_{k\el}(x).\tag{2.7}$$ By (2.6) and (2.7) we have that
$\fM^{\dt}_{\vr k\el} \fa(x)\ls 2^{-k(\f{n-1}{2 p'})}
P_{k\el}(x)$.

Finally, let $\fa$ be a $(p,n(1/p'-1))$-atom which is supported in
that ball $B(x_0;s)$. Without loss of generality, we assume that
$x_0=0$. Let $\fb(x)=s^{n/p}\,\fa(A_s\,x)$. Then $\fb$ is clearly
a $(p,n(1/p'-1))$-atom supported in the unit ball $B(0;1)$. We
also observe that $$ \split \cH^{\dt,1/t}_{\vr k\el}
*\fa(x)&=\int_{\BR^n}\cH^{\dt}_{\vr
k\el}(A_{1/t}\,x-y)\,\fa(A_t\,y)\,dy \\
&=s^{-n/p}\int_{\BR^n}\cH^{\dt}_{\vr k\el}(A_{s/t}
A_{1/s}\,x-y)\,\fb(A_{t/s}\,y)\,dy \\
&=s^{-n/p}(t/s)^{-n}\int_{\BR^n}\cH^{\dt}_{\vr
k\el}(A_{s/t}(A_{1/s}\,x-y))\,\fb(y)\,dy \\
&=s^{-n/p}\cH^{\dt,s/t}_{\vr k\el}*\fb(A_{1/s}\,x).\endsplit
\tag{2.8}$$ Therefore, combining this with the above estimate, we
complete the part (a).

(b) We observe that there is a constant $C=C(n,p)>0$ such that
for any $x_0\in\BR^n$ and for any $k\in\BN$, $\el=1,2,\cdots,N_k$,
$$\lt\|P_{k\el}(\cdot-x_0)\rt\|_{L^p}\le C. \tag{2.9}$$ Then it
easily follows from the change of variable and (2.9) that
$$\ds\lt\|(\fM^{\dt}_{\vr
k\el}\fa)\chi_{{B(x_0;2s)}^c}\rt\|_{L^p}\le C\,2^{-k(\f{n-1}{2
p'})}\|P_{k\el}(\cdot-x_0/s)\|_{L^p}\le C\,2^{-k(\f{n-1}{2
p'})}.\,\,\qed$$

{\bf Proof of Theorem 1.2.} First of all, we prove that if
$\dt>n(1/p-1/2)-1/2$ for $0<p<1$ then $\fM^{\dt}_{\vr}\fa\in
L^p(\BR^n)$ for any $(p,n(1/p'-1))$-atom on $\BR^n$ where $p'<p$
is a positive number satisfying $\dt=n(1/p'-1/2)-1/2$, and
moreover there is a constant $C>0$ independent of such atoms such
that $\|\fM^{\dt}_{\vr}\fa\,\|_{L^p}\le C$. For $t>0$ and
$\dt>0$, let
$\cH^{\dt}_{\vr,t}(x)=\cF^{-1}[(1-\vr/t)^{\dt}_+](x)$ and let
$\cH^{\dt}_{\vr,1}(x)=\cH^{\dt}_{\vr}(x)$. Let $\fa$ be a
$(p,n(1/p'-1))$-atom supported in the ball $B(x_0;s)$ with center
$x_0\in\BR^n$ and radius $s>0$. Then we see that
$\fR^{\dt}_{\vr,t}\fa(x)=\cH^{\dt}_{\vr,t}* \fa(x)$. Since
$\cH^{\dt}_{\vr}\in L^1(\BR^n)$ by Lemma 2.2, if $x\in B(0;2s)$
is given then we have that $$\lt|\fR^{\dt}_{\vr,t} \fa(x)\rt|\le
\lt\|\cH^{\dt}_{\vr,t}\rt\|_{L^1}\,\|\fa\|_{L^{\iy}}\le
\lt\|\cH^{\dt}_{\vr}\rt\|_{L^1}\, \lt|B(x_0;s)\rt|^{-1/p},$$ and
so $$\fM^{\dt}_{\vr}\fa(x)\ls \lt|B(x_0;s)\rt|^{-1/p}.$$ Since
$0<p<1$, it easily follows from (b) of Lemma 2.3 that $$ \split
\|\fM^{\dt}_{\vr}\fa\,\|^p_{L^p}&=\|(\fM^{\dt}_{\vr}\fa)\chi_{B(x_0;2s)}\|^p_{L^p}
+\|(\fM^{\dt}_{\vr}\fa)\chi_{{B(x_0;2s)}^c}\|^p_{L^p} \\
&\le 2^n +\sum_{k=1}^{\iy}\sum_{\el=1}^{N_k} \|(\fM^{\dt}_{\vr
k\el}\fa)\chi_{{B(x_0;2s)}^c}\|^p_{L^p} \\
&\ls 2^n +C \,\sum_{k=1}^{\iy} 2^{-k(\f{p}{p'}-1)(\f{n-1}{2})}\le
C.
\endsplit \tag{2.10}$$ Finally, if $f=\sum_{j=1}^{\iy} c_j \fa_j$ where
the $\fa_j$'s are $(p,n(1/p'-1))$-atoms and $\{c_j\}\in\el^p$,
then by (2.10) we have the estimate $$\|\fM^{\dt}_{\vr}
f\,\|^p_{L^p}\le \sum_j
|c_j|^p\,\|\fM^{\dt}_{\vr}\fa_j\,\|^p_{L^p}\ls \sum_j |c_j|^p.$$
Hence this completes the proof. \qed

\subheading{3. $(H^p,L^{p,\iy})$-estimate for the case that
$\Sm_{\vr}$ is a smooth convex hypersurface of finite type}

In this section we shall focus upon obtaining
$(H^p,L^{p,\iy})$-mapping properties of the maximal operator
$\fM^{\dt(p)}_{\vr}$, $p<1$, under the condition that $\Sm_{\vr}$
is a smooth convex hypersurface of finite type.

Let $\Sm$ be a smooth convex hypersurface of $\BR^n$ and let
$d\sm$ be the induced surface area measure on $\Sm$. Let
$\cE(\Sm)$ be the set of points of $\Sm$ at which the Gaussian
curvature $\kp$ vanishes, and let $\cN(\Sm)=\{n(\xi)|\,\xi\in
\cE(\Sm)\}$ where $n(\xi)$ denotes the outer unit normal to $\Sm$
at $\xi\in\Sm$. For $x\in\BR^n$, denote by $d(x/|x|,\cN(\Sm))$ the
geodesic distance on $S^{n-1}$ between $x/|x|$ and $\cN(\Sm)$, and
by $\cB (\xi(x),s)$ the spherical cap near $\xi(x)\in\Sm$ cut off
from $\Sm$ by a plane parallel to $T_{\xi(x)}(\Sm)$ ( the affine
tangent plane to $\Sm$ at $\xi(x)$ ) at distance $s>0$ from it;
that is, $$\cB
(\xi(x),s)=\{\xi\in\Sm|\,d(\xi,T_{\xi(x)}(\Sm))<s\},$$ where
$\xi(x)$ is the point of $\Sm$ whose outer unit normal is in the
direction $x$. These spherical caps play an important role in
furnishing the decay of the Fourier transform of the measure
$d\sm$. It is well known [7,9] that the function
$$\Om(\th)\fd\sup_{r>0}\,\sm[\cB(\xi(r\th),1/r)](1+r)^{\f{n-1}{2}}\tag{3.1}$$
is bounded on $S^{n-1}$ provided that $\Sm$ has nonvanishing
Gaussian curvature.

\proclaim{Definition 3.1} $\Sm$ be a smooth convex hypersurface of
$\BR^n$. Then we say that $\Sm$ satisfies a spherically integrable
condition of order $<1$ if $\Om\in L^p(S^{n-1})$ for any $p<1$.
\endproclaim

\noindent{\it Remark.} (i) B. Randol [7] proved that if $\Sm$ is a
real analytic convex hypersurface of $\BR^n$ then $\Om\in
L^p(S^{n-1})$ for some $p>2$. Thus any real analytic convex
hypersurface satifies a spherically integrable condition of order $<1$.

(ii) Let $\Sm$ be a smooth convex hypersurface of finite type
$k\ge 2$ and suppose that $\cN(\Sm)$ is a $m$-dimensional
submanifold of $\BR^n$ which is on $S^{n-1}$, where
$m<[k(n-1)]/[2(k-1)]$. Then we see ( refer to [4] ) that $\Sm$
satisfies a spherically integrable condition of order $<1$.
Moreover, it is not hard to see that $\Sm$ satisfies a spherically
integrable condition even for $m\le n-2$. We mention for reader
that it can be shown by Lemma 2.8 [4] and the fact $\Sm$ is of
finite type $P(k)$; i.e. there is some constant $C=C(\Sm)>0$ such
that for any $\th\in S^{n-1}$,
$$\Om(\th)\le\f{C}{{d(\th,\cN(\Sm))}^{\f{k-2}{2(k-1)} (n-1)}}.$$
Since $\Sm$ is smooth and of finite type, it is absolutely
impossible that $\cN(\Sm)$ is a $(n-1)$-dimensional submanifold
of $\BR^n$ which is on $S^{n-1}$.

(iii) More generally, it was shown by I. Svensson [13] that if
$\Sm$ is a smooth convex hypersurface of finite type $k\ge 2$
then $\Om\in L^p(S^{n-1})$ for some $p>2$.

Thus, by the above remark (iii), it is natural for us to obtain
the following lemma.

\proclaim{Lemma 3.2} Any smooth convex hypersurface of finite type
always satisfies a spherically integrable condition of order $<1$.
\endproclaim

Sharp decay estimates for the Fourier transform of surface measure
on a smooth convex hypersurface $\Sm$ of finite type $k\ge 2$ has
been obtained by Bruna, Nagel, and Wainger [1]; precisely
speaking, $|\cF[d\sm](x)|$ is equivalent to
$\sm[\cB(\xi(x),1/|x|)]$. They define a family of anisotropic
balls on $\Sm$ by letting $$\cB
(\xi_0,s)=\{\xi\in\Sm|\,d(\xi,T_{\xi_0}(\Sm))<s\}$$ where
$\xi_0\in \Sm$. We now recall some properties of the anisotropic
balls $\cB(\xi_0,s)$ associated with $\Sm$. The proof of the
doubling property in [1] makes it possible to obtain the following
stronger estimate for the surface measure of these balls; $$
\sm[\cB(\xi_0,\gm \,s)]\ls \lt\{ \aligned \gm^{\f{n-1}{2}}
\sm[\cB(\xi_0,s)],\,\,\,\,&\gm\ge 1, \\ \gm^{\f{n-1}{k}}
\sm[\cB(\xi_0,s)],\,\,\,\,&\gm< 1. \endaligned \rt. \tag{3.2} $$
It also follows from the triangle inequality and the doubling
property [1] that there is a positive constant $C>0$ independent
of $s>0$ such that $$\f{1}{C}\,\sm[\cB(\xi_0,s)]\le
\sm[\cB(\xi,s)]\le C\,\sm[\cB(\xi_0,s)]\,\,\,\text { for any
$\xi\in\cB(\xi_0,s)$.} \tag{3.3}$$

Next we recall a useful lemma [10] due to E. M. Stein, M. H.
Taibleson, and G. Weiss on summing up weak type functions.

\proclaim{Lemma 3.3} Let $0<p<1$. Suppose that $\{\fh_k\}$ is a
sequence of measurable functions such that for all $k\in\BN$,
$$\|\fh_k\|_{L^{p,\iy}}\le 1. $$ If $\{c_k\}\in\el^p$, then we
have the following estimate $$\lt\|\sum_{k=1}^{\iy}
c_k\fh_k\rt\|_{L^{p,\iy}}\le \lt(\f{2-p}{1-p}\rt)^{1/p}
\|\{c_k\}\|_{\el^p}.$$
\endproclaim

We now state an elementary lemma without proof which will be
useful to measure the distance from a point of $\cB(\xi_0,s)$ to
the affine tangent plane to $\Sm$ at $\xi_0\in\Sm$ in higher
dimensions.

\proclaim{Lemma 3.4} Let $\Sm$ be a smooth simple closed convex
curve in $\BR^2$ whose graph near $(0,0)$ is given as $(t,g(t))$
where $g(t)=b\,t^m +c$ is a convex function defined on $[-d,d]$
for some sufficiently small constant $b, c, d>0$ and an integer
$m\ge 2$. For $|t|\le d$, we denote by $\Th(t)$ the angle between
$n(0,g(0))$ and $n(t,g(t))$. For some small angle $\Th_0>0$ with
$\Th_0\le \max\{\Th(-d),\Th(d)\}$, let $t_0$ be chosen so that
$\Th(t_0)=\Th_0$ and $|t_0|\le d$. Then we have the following
estimate $$\lt |g(t_0)-c\rt|\sim |b|^{-\f{1}{m-1}}
m^{-\f{m}{m-1}} \Th_0^{\f{m}{m-1}}.$$
\endproclaim

\proclaim{Lemma 3.5} Let $\Sm$ be a smooth convex hypersurface of
$\BR^n$ which is of finite type $k\ge 2$. Then there is a
constant $C=C(\Sm)>0$ such that for any $y\in B(0;s)$ and $x\in
{B(0;2s)}^c$, $0<s\le 1$, $$\xi(x-y)\in\cB(\xi(x),C/|x|)$$ where
$\xi(x)$ is the point of $\Sm$ whose outer unit normal is in the
direction $x$.
\endproclaim

\noindent{\it Proof.} We observe that the following inequality
always holds for any $x, y\in \BR^n$ with $|x|>2|y|$;
$$\lt|\f{x-y}{|x-y|}-\f{x}{|x|}\rt|\le
2\,\f{|y|}{|x|}.\tag{3.4}$$ Near $\xi(x/|x|)\in \Sm$, the
hypersurface $\Sm$ can be given as the graph of a smooth convex
function defined on $D\fd T_{\xi(x/|x|)}(\Sm)\cap
B(\xi(x/|x|);1/2)$; to be precise, let $\Psi$ be a smooth convex
function defined on $D$ such that
$(\xi'_0,\Psi(\xi'_0))=\xi(x/|x|)$, and for $|t|<1/2$ and $\e\in
T^{n-1}\fd [T_{\xi(x/|x|)}(\Sm)-\xi(x/|x|)]\cap S^{n-1}$,
$$\Psi(\xi'_0+t\e)=\sum_{i=0}^k \ds\f{1}{i!}\,
\cD^i_{\e}\Psi(\xi'_0)\, t^i+\cO(t^{k+1}).\tag{3.5}$$
Using (3.5), we now estimate the distance from
$\xi=(\xi',\Psi(\xi'))\in\Sm$ to the tangent space
$T_{\xi(x/|x|)}(\Sm)$ as follows; since $\Sm$ is of finite type
$k\ge 2$, for each $\e\in T^{n-1}$ there is an integer $m$ with
$2\le m\le k$ such that for $-1/2<t<1/2$
$$\Psi(\xi'_0+t\e)-\Psi(\xi'_0)-\cD_{\e}\Psi(\xi'_0)\,t
=\ds\f{1}{m!}\,\cD_{\e}^m\Psi(\xi'_0)\,t^m +\cO(t^{m+1}).$$ Thus
by (3.4) and Lemma 3.4 we have that $$ \split \lt\la
\ds\xi\lt(\f{x}{|x|}\rt)-\xi\lt(\f{x-y}{|x-y|}\rt),\f{x}{|x|}
\rt\ra &=
\Psi(\xi'_0+t_1\e)-\Psi(\xi'_0)-\cD_{\e}\Psi(\xi'_0)\,t_1 \\
&\ls \lt[\f{m^m}{m!}\,|\cD_{\e}^m\Psi(\xi'_0)|\rt]^{-\f{1}{m-1}}
\,\lt|\f{x-y}{|x-y|}-\f{x}{|x|}\rt|^{\f{m}{m-1}} \\
&\le M_0\,\lt|\f{x-y}{|x-y|}-\f{x}{|x|}\rt| \le \f{2 M_0}{|x|}
\endsplit $$ where $t_1,\, |t_1|<1/2$, is some number so that
$(\xi'_0+t_1\,\e,\Psi(\xi'_0+t_1\,\e))=\xi\lt(\ds\f{x-y}{|x-y|}\rt)$
and $M_0=\ds\sup_{2\le m\le k}\sup_{\e\in T^{n-1}}
\lt[\f{m^m}{m!}\,|\cD_{\e}^m\Psi(\xi'_0)|\rt]^{-\f{1}{m-1}}$.
Hence we complete the proof. \qed

\proclaim{Lemma 3.6} Let $\Sm$ be a smooth convex hypersurface of
$\BR^n$ which is of finite type $k\ge 2$. Then there is a
constant $C=C(\Sm)>0$ such that for any $y\in B(0;s)$ and $x\in
{B(0;2s)}^c$, $0<s\le 1$, $$\Om\lt(\f{x-y}{|x-y|}\rt)\le
C\,\Om\lt(\f{x}{|x|}\rt)$$ where $\Om$ is the radial function
defined as in (3.1).
\endproclaim

\noindent{\it Proof.} It easily follows from (3.2), (3.3), the
definition of $\Om$, and Lemma 3.5 that for any $y\in B(0;s)$ and
$x\in {B(0;2s)}^c$, $0<s\le 1$, $$ \split
\Om\lt(\f{x-y}{|x-y|}\rt)&=\sup_{r>0}\sm[\cB(\xi(x-y),1/r)]\,
(1+r)^{\f{n-1}{2}} \\
&\ls \sup_{r>0}\sm[\cB(\xi(x),1/r)]\,
(1+r)^{\f{n-1}{2}}=\Om\lt(\f{x}{|x|}\rt).\,\,\, \qed \endsplit $$

{\bf Proof of Theorem 1.1.} Fix $0<p<1$. Let $\fa$ be a
$(p,n(1/p-1))$-atom supported in the ball $B(x_0;s)$ with center
$x_0\in\BR^n$ and radius $s>0$. Then we see that
$\fR^{\dt}_{\vr,t}\fa(x)=\cH^{\dt}_{\vr,t}* \fa(x)$. Recalling
the lemma [6] about asymptotics of quasiradial Bochner-Riesz
kernel and the result of Bruna, Nagel, and Wainger [1], we get
that $$\lt|\cH^{\dt(p)}_{\vr}(x)\rt|\sim
\lt|\n\cH^{\dt(p)}_{\vr}(x)\rt|
\sim\f{1}{(1+|x|)^{\f{n}{p}-\f{n-1}{2}}} \,\sm[\cB(\xi(x),1/|x|)]
\tag{3.6}$$ where we consider $\Sm_{\vr}$ as $\Sm$ given in the
above. Since $\cH^{\dt(p)}_{\vr}\in L^1(\BR^n)$ by
 (3.6) and Lemma 3.2, if $x\in B(0;2s)$ is given then we have that
$$\lt|\fR^{\dt(p)}_{\vr,t} \fa(x)\rt|\le
\lt\|\cH^{\dt(p)}_{\vr,t}\rt\|_{L^1}\,\|\fa\|_{L^{\iy}}\le
\lt\|\cH^{\dt(p)}_{\vr}\rt\|_{L^1}\, \lt|B(x_0;s)\rt|^{-1/p},$$
and so $$\fM^{\dt(p)}_{\vr}\fa(x)\ls \lt|B(x_0;s)\rt|^{-1/p}.$$
Thus we have that for all $\ld>0$, $$\lt|\{x\in
B(x_0;2s)|\,\fM^{\dt(p)}_{\vr}\fa(x)>\ld/2\}\rt|\ls
\ld^{-p}.\tag{3.7}$$

Next we shall obtain the following inequality $$\lt|\{x\in
{B(x_0;2s)}^c|\,\fM^{\dt(p)}_{\vr}\fa(x)>\ld/2\}\rt|\ls
\ld^{-p},\,\,\ld>0.\tag{3.8}$$

As in the argument of (2.8), without loss of generality we can
assume that a $(p,n(1/p-1))$-atom $\fa$ is supported in the unit
ball $B(0;1)$ centered at the origin. We now consider the case
that $x\in {B(0;2)}^c$ and $t>1$. Then it follows from (3.1),
(3.2), (3.6), and Lemma 3.6 that $$ \split
\lt|\cH^{\dt(p)}_{\vr,t} *\fa(x)\rt| &\ls t^n\int_{B(0;1)}
|\cH^{\dt(p)}_{\vr}(A_t(x-y))|\,dy \\
&\ls \f{t^{n-n/p}}{(1+|x|)^{\f{n}{p}}}
\int_{B(0;1)}\Om\lt(\f{x-y}{|x-y|}\rt)\,dy
\\ & \ls \f{t^{n-n/p}}{(1+|x|)^{\f{n}{p}}}
\,\Om\lt(\f{x}{|x|}\rt) \\ &\ls \f{1}{(1+|x|)^{\f{n}{p}}}
\,\Om\lt(\f{x}{|x|}\rt)
\endsplit $$ because $n(1-1/p)<0$. So we have that
$$\sup_{t>1}\lt|\cH^{\dt(p)}_{\vr,t} *\fa(x)\rt|\ls
\f{1}{(1+|x|)^{\f{n}{p}}}\,\Om\lt(\f{x}{|x|}\rt). \tag{3.9} $$

Let $N\in\BN$ be an integer satisfying $N<n(1/p-1)\le N+1$, i.e.
$n/(n+N+1)\le p<n/(n+N)$. If $x\in {B(0;2)}^c$ and $0<t\le 1$, let
$\cQ_{t,x}(y)$ be the $N$-th order Taylor polynomial of the
function $y\mapsto \cH^{\dt(p)}_{\vr}(A_t^*(x-y))$ expanded near
the origin, where
$\cH^{\dt(p)}_{\vr}(x)=\cF^{-1}[(1-\vr)^{\dt(p)}_+](x)$. Then it
follows from the moment condition on the atom $\fa$, Taylor's
theorem, (3.1), (3.2), (3.6), and Lemma 3.6 that $$ \split
\lt|\cH^{\dt(p)}_{\vr,t} *\fa(x)\rt| &=
t^n\lt|\int_{\BR^n}[\cH^{\dt(p)}_{\vr}(A_t(x-y))-\cQ_{t,x}(y)]\fa(y)\,dy
\rt| \\ &\ls t^{n+(N+1)}\int_0^1\int_{B(0;1)}
|\n^{N+1}\cH^{\dt(p)}_{\vr}(A_t(x-\tau y))|\,dy\,d\tau \\
&\ls \f{t^{n+(N+1)-n/p}}{(1+|x|)^{\f{n}{p}}}
\int_0^1\int_{B(0;1)}\Om\lt(\f{x-\tau y}{|x-\tau
y|}\rt)\,dy\,d\tau
\\ & \ls \f{t^{n+(N+1)-n/p}}{(1+|x|)^{\f{n}{p}}}
\,\Om\lt(\f{x}{|x|}\rt) \\ &\ls \f{1}{(1+|x|)^{\f{n}{p}}}
\,\Om\lt(\f{x}{|x|}\rt)
\endsplit $$ because $n+(N+1)-n/p\ge 0$. Thus we have that
$$\sup_{0<t\le 1}\lt|\cH^{\dt(p)}_{\vr,t} *\fa(x)\rt|\ls
\f{1}{(1+|x|)^{\f{n}{p}}}\,\Om\lt(\f{x}{|x|}\rt). \tag{3.10}$$
Thus by (3.9) and (3.10) we conclude that $$\fM^{\dt(p)}_{\vr}
\fa(r\th)\ls \f{1}{(1+r)^{\f{n}{p}}}\,\Om(\th).$$ Hence we have
the following estimate $$\int_{\{x\in
{B(0;2)}^c|\,\fM^{\dt(p)}_{\vr} \fa(x)>\ld\}}\,dx \ls
\int_{S^{n-1}}\int_{\{r>0|\,2<r<\ld^{-p/n}\,
{\Om(\th)}^{p/n}\}}\,r^{n-1}\,dr\,d\th \ls \ld^{-p}$$ because
$\Om\in L^p(S^{n-1})$ for any $p<1$ by Lemma 3.2. Therefore, by (3.7), (3.8),
and Lemma 3.3, we complete the proof. \qed

\noindent{\it Acknowledgements.} The author had a chance to
present this manuscript in Workshop on Fourier Analysis and
Convexity, Universit\`a di Milano-Bicocca, Italy, June 11-22,
2001, organized by Professors Brandolini, Colzani, Iosevich, and
Travaglini. He would like to thank for their hospitality and
kindness during having stayed there, and also had a wonderful
impression for friendship of lots of participants from all over
the world. Especially, it was a great pleasure to have a chance to
discuss with Professor Terry Tao on the unsolved problem which he
mentioned in Theorem 1.2. The author would not have got a clue
without stimulating discussion with him, and would like to thank
for his kindness and concern. Finally the author would like to
thank Professor Galia Dafni for her concern on this subject.

\enddocument